\documentclass [11pt]{article}		

\usepackage{graphicx}
\usepackage{amsmath}
\usepackage{amssymb}

\newcommand{\bdry}{\partial}

\renewcommand{\Hat}{\widehat}

\newcommand{\pf}[1]{\noindent{\bf Proof#1:\  }}
\newcommand{\defn}[1]{\noindent{\bf Definition#1:\  }}

\newcommand{\eg}[1]{\noindent{\bf Example#1:\  }}

\newcommand{\QED}{\hfill$\square$\medskip}

\newtheorem{thm}{Theorem}[section]
\newtheorem{prop}[thm]{Proposition}
\newtheorem{lemma}[thm]{Lemma}

\begin{document}

\title{A new partition identity coming from complex dynamics}
\author{George E. Andrews\footnote{Department of Mathematics \& Statistics, Penn State University, University Park, State College, PA, 16802. USA.  andrews@math.psu.edu Research supported in part by a grant from NSF.} \and Rodrigo A. P\'erez\footnote{Department of Mathematics, Cornell University, Ithaca, NY 14853. USA.  rperez@math.cornell.edu Partially supported by NSF Postdoctoral Fellowship, grant DMS-0202519.}}
\date{}

\maketitle

\begin{abstract}
{\it We present a new identity involving compositions (i.e. ordered partitions
  of natural numbers). The Formula has its origin in complex dynamical systems
  and appears when counting, in the polynomial family $\{ f_c:z \mapsto z^d +
  c \}$, periodic critical orbits with equivalent itineraries. We give two
  different proofs of the identity; one following the original approach in
  dynamics and another with purely combinatorial methods.}
\end{abstract}

\section{Introduction}\label{sect:Intro}
The field of dynamical systems takes frequent advantage of combinatorial
techniques to classify all sorts of dynamic phenomena. Often the tools
borrowed are classic, so there are few opportunities for feedback. In this
note, we introduce a previously unknown identity in the theory of partitions,
which arose from dynamical considerations. We will give two different proofs
of the formula; one that illustrates the original approach in dynamics and
a second one using the more traditional methods of enumerative
combinatorics. \\

\defn{s} Let $n \in \mathbb{N}$. A {\it composition} of $n$ in $r$ parts is a
  partition $P:=[a_1 + \ldots + a_r = n]$ that takes into account the order of
  the parts $a_j.$ An {\it H-composition} is a composition satisfying $a_1
  \geq a_j$ for all $j \leq r$. We use $\mathcal{H}(n)$ to denote the
  collection of H-compositions of $n.$

  The {\it multiplicity} $\omega$ of $P \in \mathcal{H}(n)$ is defined as the
  number of parts other than $a_1$, equal to $a_1;$ that is, $\omega(P) = \#
  \{ j > 1 \mid a_j = a_1 \}$. \\

\renewcommand{\theequation}{{\bf \arabic{section}.\arabic{equation}}}
{
\begin{thm}\label{Thm:Formula}
  For all $n,d \in \mathbb{N}$ the following identity holds
  {
  \begin{equation}\label{eqn:Formula}
    \sum_{P \in \mathcal{H}(n)}
      \varphi(a_1) \cdot (d-1)^{r-\omega(P)} \cdot d^{\omega(P)} =
    d^n-1
  \end{equation}
  }
  where $\varphi$ represents Euler's totient function.
\end{thm}
}

\eg{}
  $\mathcal{H}(5) = \{ 5, 4+1, 3+2, 3+1+1, 2+2+1, 2+1+2, 2+1+1+1, 1+1+1+1+1
  \}$. For degree $d=3$, the formula yields:
  \begin{gather*}
    (4 \cdot 2 \cdot 1) + (2 \cdot 2^2 \cdot 1) + (2 \cdot 2^2 \cdot 1) + 
    (2 \cdot 2^3 \cdot 1) + (1 \cdot 2^2 \cdot 3) + \\
    + (1 \cdot 2^2 \cdot 3) + (1 \cdot 2^4 \cdot 1) + (1 \cdot 2 \cdot 3^4) =
    242 = 3^5 - 1
  \end{gather*}

Formula \ref{eqn:Formula} was first detected in an effort to list the possible
combinatorial behaviors of critically periodic orbits in the family of complex
polynomials $\{ f_c:z \mapsto z^d + c \mid c \in \mathbb{C}\}$ of fixed degree
$d$. Every polynomial function $f:\mathbb{C} \longrightarrow \mathbb{C}$ has
associated a compact set $K$, its filled Julia set, that is invariant under
$f$. When the critical point is periodic, $f$ is described by a finite amount
of data that encodes the location in $K$ of points in the orbit of 0. Theorem
\ref{Thm:Formula} is proved by counting polynomials with equivalent
descriptions.

In Section \ref{sec:basics} we provide a condensed review of the relevant
concepts from complex dynamics. This will furnish a language to describe the
dynamical picture and give some intuition on the behavior of critical
orbits. Admittedly, the statements that we need to quote far exceed our
limitations of space, and the consequence is a constant referral to the
literature. We would like to call attention to \cite{Eber} and
\cite{Poir}. These works deserve more publicity as they clarify the status of
many folk results that had no prior reference.

Section \ref{sec:count} uses the material introduced before to prove Theorem
\ref{Thm:Formula} from the viewpoint of complex dynamics. Even though the
proof of several supporting claims is deferred to the references, the
inclusion of this method is justified by its potential to uncover similar
identities. This is briefly mentioned at the end of the Section, where a few
remarks are made on the combinatorial structure of Formula
(\ref{eqn:Formula}).

A self-contained proof of the Formula, relying exclusively on enumerative
combinatorics, is presented in Section \ref{sec:pf2}.

\section{Basics in Complex Dynamics}\label{sec:basics}
In this section we sketch the basic material in dynamics of polynomials in one
complex variable. Proofs of the results stated and further information can be
found in \cite{Orsay}, \cite{M_book} and \cite{Eber}. The focus here will be
on binomials of the form $f_c(z)=z^d+c.$ This family covers all affine
conjugacy classes of polynomials with exactly one critical point. For any
point $z$, the sequence $\big\{ z, f_c(z), f_c^2(z), \ldots \big\}$ of $z$ is
called the {\it orbit} of $f_c$ and is denoted $\mathcal{O}(z)$.

The {\it filled Julia set} associated to $f_c$ is
\begin{equation}\label{eqn:K}
  K_c = \big\{ z \in \mathbb{C} \mid
               \mathcal{O}(z) \text{ is bounded} \big\}.
\end{equation}
$K_c$ is a perfect set; i.e. it contains all its accumulation points. It is
totally invariant under $f_c$; that is, $f_c(K_c) = f_c^{-1}(K_c) = K_c.$
Depending on whether the critical point belongs to the filled Julia set or
not, $K_c$ is simply connected or a Cantor set. Moreover, since $f_c$ is a $d$
to 1 cover of $\mathbb{C}$ branched only at the critical point, $K_c$ has
$d$-fold rotational symmetry around 0.

A point $z$ is called {\it periodic} if $f_c^n(z)=z$ for some $n \geq 1$. The
least such $n$ is called the {\it period} of $z$ and the value $\lambda =
\left(f_c^n \right)'(z)$ associated to $\mathcal{O}(z)$ is the {\it
multiplier} of the orbit. When $n=1,$ $z$ is a {\it fixed} point. A periodic
orbit is called {\it attracting, indifferent} or {\it repelling} depending on
whether $|\lambda|$ is less than, equal or greater than 1. Note that when the
critical point belongs to a periodic orbit, the multiplier is 0; we speak then
of an {\it superattracting} orbit or say that the map $f_c$ has {\it periodic
critical orbit.}


With the exception of $z^2 -\frac34,$ every binomial $z^d+c$ has at least one
periodic orbit of every period\footnote{For quadratic binomials there is a
unique orbit of period 2. As $c\rightarrow -\frac34,$ both points in the
period 2 orbit of $z^2+c$ approach each other and collapse into a fixed point
of multiplier -1. In all other cases, even if one orbit collapses as the
parameter varies, there are other orbits of the same period that persist.}. In
particular, by the Fundamental Theorem of Algebra, $f_c$ always has $d$ fixed
points counted with multiplicity.

Most elementary dynamical properties can be deduced from the behavior of
critical points and their relation to periodic orbits. The pioneering work of
P. Fatou and G. Julia around 1918 produced the following results (valid for
arbitrary complex polynomials; see \cite{Orsay} or \cite{M_book}):
{\renewcommand{\theenumi}{{\bf J\arabic{enumi}}}
\begin{enumerate}
  \item \label{J-one} To every attracting orbit $\mathcal{O}(z)$ corresponds
        at least one critical point $c_0$ such that $\mathcal{O}(c_0)$ is
        captured in a neighborhood of $\mathcal{O}(z)$ and eventually
        converges to this orbit.
  \item \label{J-two} Attracting orbits are contained in the interior of the
        filled Julia set $K,$ whereas repelling orbits belong to $\bdry K.$
        Moreover, $\bdry K$ is the closure of the union of all repelling
        orbits.
\end{enumerate}
}

In the case that we study, the only critical point of $f_c(z)=z^d+c$ is 0, so
statement \ref{J-one} implies that $f_c$ can have at most one attracting orbit
$\mathcal{O}_{\text{attr}}$. If this orbit exists, the iterates of 0 converge
to $\mathcal{O}_{\text{attr}}$; then, $\mathcal{O}(0)$ is bounded and it
follows from (\ref{eqn:K}) that $K_c$ is simply connected.

\begin{figure}[h]
  \hspace{1.3cm}
  \includegraphics{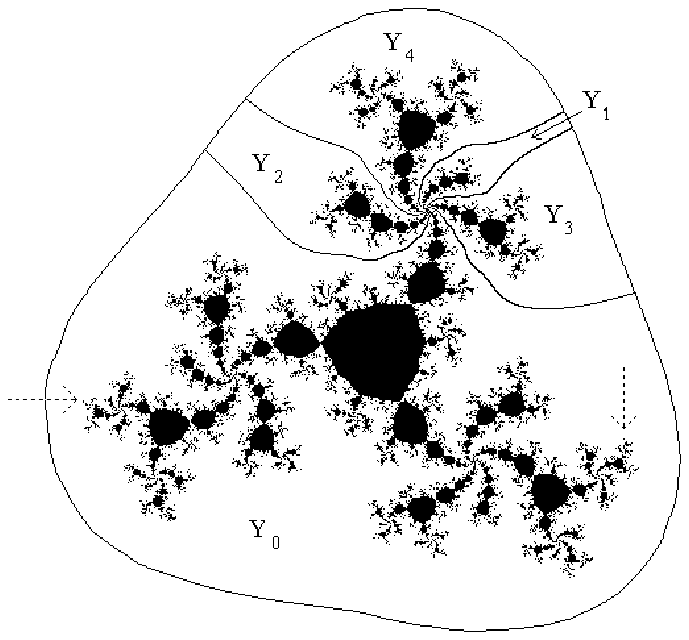}
  \hspace{.3in}
  \includegraphics{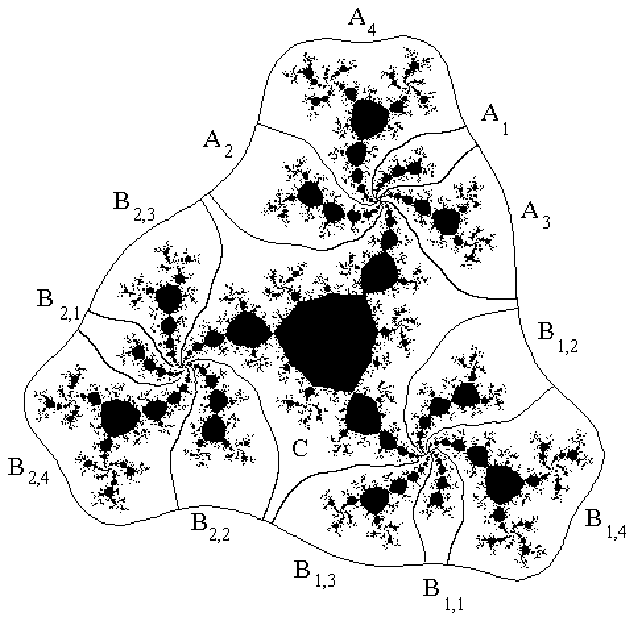}
  \caption{\small Let $c=(.387848...)+i(.6853...).$ The left picture shows the
    filled Julia set $K_c$ of the cubic map $z^3 + c,$ covered by level 0 of
    the puzzle. The center of symmetry is at 0, the point where the rays
    converge is $\alpha$ and the other fixed points are marked by dotted
    arrows. In this example the rotation number around $\alpha$ is
    $\rho_{\alpha} = \frac25$ and the ray angles are $\frac{5}{121} \mapsto
    \frac{15}{121} \mapsto \frac{45}{121} \mapsto \frac{14}{121} \mapsto
    \frac{42}{121} \mapsto \frac{5}{121}$. The right picture illustrates level
    1 of the puzzle for the same map.}
  \label{fig:Puzzles}
\end{figure}

In the remainder of this Section $c$ is chosen so that
$\mathcal{O}_{\text{attr}}$ exists and has period at least 2. Then all fixed
points will be repelling and in particular, belong to $\bdry K_c$. To
understand better the geometry of $K_c$ consider the left picture in Figure
\ref{fig:Puzzles}. Of the $d$ fixed points, exactly one, denoted $\alpha,$ has
the property that $K_c \setminus \alpha$ splits into several disjoint
components. The component that contains 0 will contain also the remaining
$d-1$ fixed points.

Notice that $K_c$ has a ''fractal structure''. This is illustrated for
instance by the fact that all $n$-fold preimages of $\alpha$ separate $K_c$
into as many components as $\alpha$ does. Moreover, such preimages are dense
in $\bdry K_c$.

Let $\phi_c:\Hat{\mathbb{C}} \setminus K_c \longrightarrow \Hat{\mathbb{C}}
\setminus \mathbb{D}$ be the Riemann map between the complements of $K_c$ and
the unit disk, normalized to have derivative 1 at $\infty.$ The pull-back by
$\phi_c$ of concentric circles ($|\zeta| = r,$ with $\zeta \in \mathbb{C}$)
yields a family of {\it equipotential curves} enclosing $K_c.$ Similarly, the
pull-back of radial lines ($\arg \zeta = \theta$) results in a family of {\it
exterior rays} emanating from $K_c.$ These two families of curves form
mutually orthogonal foliations of $\mathbb{C} \setminus K_c.$ The
equipotential curve of radius $r>1$ and the external ray of angle $\theta$
will be denoted by $e_r$ and $r_{\theta}$ respectively.

The appeal of working with foliations by equipotentials and external rays lies
in the fact that they are invariant under the action of $f_c;$ more precisely,
we have the relations
\begin{equation}\label{eqn:Foliations}
  f_c(e_r) = e_{(r^d)}
  \hspace{.3in}
  f_c(r_{\theta}) = r_{(d\theta)}
\end{equation}

It is important to point out that the normalization of $\phi_c$ determines the
branch of $f_c^{-1}(r_{(d\theta)})$ that corresponds to $r_{\theta}$. Property
(\ref{eqn:Foliations}) of the equipotential and ray foliations is the basis
for the definition of the {\it Yoccoz puzzle}: Fix the neighborhood $U$ of
$K_c$ bounded by the equipotential of radius 2 (any other radius $>1$ will do)
and consider the collection $\mathcal{R}_\alpha$ of rays landing\footnote{A
ray $r_{\theta}$ lands at a point $z \in K_c$ if $z$ is the only accumulation
point of $r_{\theta}$ in $K_c$. The issue of landing is a delicate one as rays
could accumulate on a large subset of $K_c$. However, rays with rational
angles always land at a unique point and this is enough for our purposes
(\cite{Orsay}, \cite{M_book}).} at $\alpha$; refer to Figure
\ref{fig:Puzzles}.

$\mathcal{R}_{\alpha}$ is a finite set and it is known that $f_c$ acts on it
by a cyclic permutation. If each ray is sent counterclockwise to the ray $p$
positions ahead, the {\it rotation number around $\alpha$} is given by
$\rho_{\alpha} = \frac pq$ where $q=|\mathcal{R}_\alpha|,$ $p<q$ and
$\text{gcd}(p,q) = 1.$ The rays in $\mathcal{R}_\alpha$ have rational angles
that depend on the values $\rho_{\alpha}$ and $d$; they split the region $U$
into $q$ connected components whose closures will be called the {\it puzzle
pieces} of {\it level} $0$ and denoted $Y_0, Y_1, \ldots, Y_{q-1}.$ Here the
subindices are residues modulo $q$ and are chosen so that $0 \in Y_0$ and $f_c
\left( K_c \cap Y_j \right) = K_c \cap Y_{j+1}$ for $j=1, \ldots, q-1$; in
particular, the critical value $c$ is in $Y_1$. The combinatorial richness
hidden in this picture follows from the fact that
\begin{equation}\label{eqn:Richness}
  f_c \left( K_c \cap Y_0 \right) = K_c,
\end{equation}
creating multiple overlappings; we expand on this situation below, where level
1 of the puzzle is discussed in more detail.

The puzzle pieces of level $m$ are defined as the closures of the connected
components of $f_c^{-m} \big( \bigcup \text{ int\,} Y_j \big)$ for $j=0\ldots
q-1.$ The resulting family $\mathcal{Y}_c$ of puzzle pieces of all levels has
the following properties:

{\renewcommand{\theenumi}{{\bf Y\arabic{enumi}}}
\begin{enumerate}
  \item \label{Y-one} Each piece is a closed topological disk whose boundary
        is formed by segments of rays landing at preimages of $\alpha$ and
        segments of an equipotential curve. To each level of the puzzle there
        corresponds one equipotential.
  \item \label{Y-two} There are $(q-1)d^m+1$ pieces of level $m$ and they form
        a covering of $K_c.$ The unique piece that contains the critical point
        is called the {\it critical piece} of level $m$.
  \item \label{Y-three} Any two puzzle pieces either are nested (with the
        piece of higher level contained in the piece of lower level), or have
        disjoint interiors.
  \item \label{Y-four} The image of any piece $Y$ of level $m \geq 1$ is a
        piece $Y'$ of level $m-1$. The restricted map $f_c:\text{int\,}Y
        \longrightarrow \text{int}Y'$ is a $d$ to 1 branched covering or a
        conformal homeomorphism, depending on whether $Y$ is critical or not.
\end{enumerate}}

Next, we will give names to the pieces of level 1 and describe briefly their
adjacencies and behavior under $f_c$; consult the right side of Figure
\ref{fig:Puzzles} for reference and \cite{Eber}, \cite{Multibrots} for
information on the case $d>2$. Let $C$ be the critical piece of level 1. $C$
has $d$-fold symmetry and the intersection of its boundary with $K_c$ consists
of all the points in $f_c^{-1}(\alpha)$, including $\alpha$ itself (see
Property \ref{Y-one}). We label these points $\alpha,\alpha_1,
\alpha_2,\ldots, \alpha_{d-1}$ as they are located clockwise on $\bdry C.$

Besides $C$, there is a fan of pieces around $\alpha$; we call them $A_1,
\ldots, A_{q-1}$. Here, labels are chosen so that $K_c \cap A_j = K_c \cap
Y_j$. Thus, $f_c(A_j) = Y_{j+1}$ ($j=1, \ldots, q-1$) is 1 to 1, while $f(C)
= Y_1$ is a $d$ to 1 branched cover; compare Property \ref{Y-four}.

The picture is similar around all $\alpha_k$. There are $q$ pieces $C,\,
B_{k,1},\, \ldots, \, B_{k,q-1}$, but here $f_c$ does not permute the pieces
around $\alpha_k$; instead $f_c(B_{k,j}) = Y_{j+1}$. In short, each $Y_j\, (j
\neq 1)$ has $d$ preimages $A_{j-1}, B_{1,j-1}, \ldots, B_{d-1,j-1}$, while
$f_c^{-1}(Y_1) = C$. Here again the indices are residues modulo $q$ so
$f_c(B_{k,q-1})=Y_0$ for any $1 \leq q \leq d-1$.

As a consequence, consider a point $z \in A_j \cap K_c$. Its orbit is forced
to cyclically follow the rest of the fan $A_{j+1}, A_{j+2}, \ldots$ around
$\alpha$ until $f^{q-j-1}(z) \in A_{q-1} \cap K_c$ and, one step later,
$f^{q-j}(z) \in Y_0 \cap K_c$. Because of (\ref{eqn:Richness}), the next
iterate could be located anywhere in $K_c$ depending on the exact position of
$f^{q-j}(z)$ within $Y_0$.

\section{Counting Hyperbolic Components}\label{sec:count}

When the critical orbit $\mathcal{O}(0)$ of $f_c$ is periodic, its behavior
can be classified according to the disposition within $K_c$ of the points in
$\mathcal{O}(0)$.

The first proof of Formula (\ref{eqn:Formula}) will be based on a careful
study of the different patterns attainable by the critical orbits of
critically periodic binomials. \\

\defn{} A {\it $d$-center} is a parameter $c$ such that the map $f_c(z)=z^d+c$
  has periodic critical orbit. We will refer to any $n$ for which $f_c^n(0)=0$
  as a {\it period} of the $d$-center.

\begin{lemma}\label{count}
  The number of different $d$-centers with period $n$ is $d^{n-1}$.
\end{lemma}

\pf{ (Gleason \cite[expos\'e XIX]{Orsay})} Given $n$ we want to count all the
  solutions $c$ of the equation $f_c^n(0) = 0$. Since $f_c(0)=c$, this is
  equivalent to count the solutions of $f_c^{n-1}(c)=0$. This is a polynomial
  of degree $d^{n-1}$ in $c$; hence we only need to show that all its
  solutions are different.

  Define the family of polynomials $\big\{ h_r \in \mathbb{Z}[z] \big\}$ by
  the recursion $h_0(z) = z$, $h_r(z) = \big( h_{r-1}(z) \big)^d + z$, so that
  the critical orbit of $f_c$ returns to 0 after $n$ iterations if and only
  ifthat our condition reads $h_{n-1}(c) = 0$. Each $h_r(z)$ is a monic
  polynomial with integer coefficients, showing that $c$ belongs to the ring
  $\mathbb{A}$ of algebraic integers.

  Suppose $c$ is a multiple root of $h_n(z)$; that is, $h_n'(c)=0$. From
  $h_n'(z) = d \big( h_{n-1}(z) \big)^{d-1} \cdot h_{n-1}'(z) + 1$ we conclude
  that $\big( h_{n-1}(c) \big)^{d-1} \cdot h_{n-1}'(c) = \frac{-1}d$. By the
  additive/multiplicative closure of $\mathbb{A}$, the left hand expression is
  again an algebraic integer; thus $\frac{-1}d \in \mathbb{Q} \cap \mathbb{A}
  = \mathbb{Z}$ ! (refer for instance to \cite{Dedekind}). This contradiction
  shows that $c$ must be a simple root of $h_n(z)$ and the result follows.
  \QED

\defn{s} Choose a $d$-center $c$ with period $n$ and let $\text{Crit}_0$ be
  the ordered set $\big\{ z_0, z_1, \ldots, z_n \big\}$ where $z_j=f_c^j(0),\,
  j=1, \ldots, n$ describes one period of the critical orbit $\mathcal{O}(0)$.
  Let $Z_j$ be the puzzle piece of level $n-j$ that contains $z_j$; in
  particular, $Z_n$ is the critical piece $Y_0$. It is convenient to think of
  the family $Z_n, \ldots, Z_0$ as defined in descending order of indices by
  the finite recursion $Z_n=Y_0,\, f_c(Z_{j-1})=Z_j\, (j=n, \ldots, 1)$.
  Accordingly we will say that $(Z_{j-1},z_{j-1})$ is the {\it pull-back} of
  $(Z_j,z_j)$ {\it along} $\text{Crit}_0$. By Property \ref{Y-four}, either
  $Z_j$ contains the critical value $c$ and $Z_{j-1}=f_c^{-1}(Z_j)$ or
  $Z_{j-1}$ is one of the $d$ pieces that constitute $f_c^{-1}(Z_j)$.

  We will associate to $c$ an {\it itinerary} $\big( (a_1,b_1), \ldots,
  (a_r,b_r) \big)$ as follows. Consider the subsequence $\zeta_0, \ldots,
  \zeta_r$ of those points in $\text{Crit}_0$ that lie in $Y_0$. In particular
  $\zeta_0=0$, $\zeta_1=z_q$ and $\zeta_r=0$ again. The numbers $a_i$ are
  defined by the relation $f_c^{a_i}(\zeta_{i-1}) = \zeta_i$. Observe that
  $a_i=q$ exactly when $\zeta_{i-1} \in C$ (this includes the case $\zeta_0=0
  \,\,\Rightarrow\,\, a_1=q$). When $a_i=q$ we let $b_i=0$. Otherwise $a_i<q$,
  and then $\zeta_{i-1} \in B_{k,q-a_i}$ for some $k$. The term $q-a_i$
  appears simply because $f_c(\zeta_{i-1}) \in f_c(B_{k,q-a_i})=A_{q-a_i+1}$
  as required by the definition of $a_i$ (compare with the discussion at the
  end of Section \ref{sec:basics}). In this situation we let $b_i=k$. Every
  pair $(a_i,b_i)$ will be called a {\it leg} of the itinerary.

  From the definition of the $a_i$ it is immediate that $f_c^n(0) = f_c^{a_r}
  \circ \ldots \circ f_c^{a_1}(0)$ and therefore $a_1+ \ldots +a_r = n$. Since
  $a_1 = q \geq a_j$ ($1<j\leq r$), it follows that $a_1+ \ldots + a_r = n$ is
  an H-composition of $n$ with $r$ parts and we denote it by $P(c)$. \\

The above definitions afford us the means to describe critical orbits. The
distribution of elements of $\text{Crit}_0$ within $K_c$ is well conveyed by
its itinerary, even though theses objects are not in 1 to 1 correspondence. The
core result in this Section makes precise just how much extra information is
required to single out a  particular $d$-center:

\begin{prop}\label{prop:Core}
  Let $P$ denote an H-composition $a_1 + \ldots + a_r = n$ with $a_1>1$ and
  multiplicity $\omega(P) = w$. Then the total number of $d$-centers $c$
  having $n$ as a period and such that $P(c) = P$ is
  \begin{equation}\label{mult_of_combs}
    \#\big\{ c\text{ \rm is a $d$-center } \mid P(c)=P \big\} =
    \varphi(a_1) \cdot (d-1)^{r-w} \cdot d^w.
  \end{equation}
\end{prop}

Theorem \ref{Thm:Formula} follows from the previous results. Lemma \ref{count}
determines the total number of $d$-centers with period $n$, while Proposition
\ref{prop:Core} sorts them by combinatorial type.

\pf{ of Theorem \ref{Thm:Formula}} Note that the binomial $f_0(z)=z^d$ is the
  only one with 0 as a fixed point. Thus, the H-composition $1+1+\ldots+1 =
  n$ can be associated to the single $d$-center $c=0$, regardless of the value
  of $n$. The other $d$-centers are classified in Proposition \ref{prop:Core}
  according to their associated H-composition, so by Lemma \ref{count} the
  total of $d$-centers is
  \begin{equation*}
    \text{\Large 1} \hspace{.2in}
    + \sum_{\substack{P \in \mathcal{H}(n) \\ a_1>1}}
    \varphi(a_1)(d-1)^{r-\omega(P)}p^{\omega(P)} = d^{n-1}.
  \end{equation*}
  Since $\omega(1+1+ \ldots +1) = n-1$, the LHS can be modified to incorporate
  this particular case under the sum symbol to coincide with the sum in
  Formula (\ref{eqn:Formula}). The adjusted value on the RHS becomes
  $d^{n-1}-1 + (d-1)d^{n-1} = d^n - 1$ as claimed.
\QED

\pf{ of Proposition \ref{prop:Core}} The proof will be divided in 2 parts
  according to the structure of the given H-composition $P$. Essentially we
  have to handle apart the possibility that $P$ admits $d$-centers with period
  smaller than $n$. Let us describe first the situation of period less than
  $n$ in order to present an outline of the proof.

  Suppose that for some $z^d+c$ and $j \leq n-1$, the piece $Z_{n-j}$ contains
  0 as well as $z_{n-j} \in \text{Crit}_0$.  Since $z_n=0$ and $Z_{n-j}
  \subsetneq Z_n$ (by Property \ref{Y-three}), it follows that $Z_{n-2j}$
  contains $z_n=0,\, z_{n-j}$ and $z_{n-2j}$. By the same argument, every
  $Z_{n-kj}$ contains $0,z_{n-j},z_{n-2j}, \ldots, z_{n-kj}$ as long as $n
  \geq kj$. More generally, the index $i:=\gcd(j,n)$ has the property that
  \begin{equation*}
    0,z_i,z_{2i}, \ldots, z_{n-2i},z_{n-i} \in Z_0 \subsetneq Z_i
    \subsetneq Z_{2i} \subsetneq \ldots \subsetneq Z_{n-i} \subsetneq Y_0
  \end{equation*}\
  and $f_c^i$ maps $Z_0 \mapsto Z_i \mapsto \ldots \mapsto Z_{n-i} \mapsto
  Y_0$. Hence, $i<n$ is a period of the $d$-center $c$. Moreover, since any 2
  points $z_{ki}, z_{(k+1)i}$ are in $Z_0$, they must follow for $i$
  consecutive steps the same itinerary. As a consequence the full itinerary
  has the following form
  \begin{equation}\label{eqn:Ren}
    \left( \overbrace{(a_1,b_1), \ldots, (a_i,b_i)\hspace{.1in},
    \hspace{.25in}.\hspace{.25in}.\hspace{.25in}.\hspace{.25in},\hspace{.1in}
    (a_1,b_1),\ldots, (a_i,b_i)}^{\text{repeated $\frac ri$ times}} \right).
  \end{equation}

  When this happens we say that the underlying H-composition is {\it
  renormalizable}. Any itinerary with the structure of (\ref{eqn:Ren}) is said
  to be {\it renormalizing}. Observe that a renormalizable H-composition may
  give rise to a non-renormalizing itinerary; it is enough that one of the
  $b_j$ does not match the pattern in (\ref{eqn:Ren}).  Additionally, for a
  renormalizing itinerary any of its associated $d$-centers has period
  $i<n$. In order to deal with these deterrents, the case of renormalizable
  H-compositions will be treated last. \\

  Our strategy is to show that every itinerary associated to the given
  H-composition $P$ corresponds to a fixed number of $d$-centers. The outline
  of the proof is as follows. If an itinerary is non-renormalizing, we count
  all pairs of angles $\eta^-, \eta^+$ such that the rays $r_{\eta^-},\,
  r_{\eta^+}$ can delimit a pull-back piece $Z_1$. By results of Goldberg and
  Milnor (\cite{Gold}, \cite{Portraits}), $d$-centers are in 1 to 1
  correspondence with such pairs of angles. If $P$ is not renormalizable,
  every itinerary is non-renormalizing and the result follows.

  In the case of renormalizable H-compositions we separate the different
  itineraries in renormalizing and non-renormalizing. Reducing every
  renormalizing itinerary to the non-renormalizing itinerary of a higher
  degree binomial, we get again the count $\varphi(a_1) \cdot (d-1)^{r-w}
  \cdot d^w$. \\

  \noindent\underline{Non-renormalizing itineraries:} Let $c$ be a $d$-center
  such that $P(c)=P$. The rotation number around $\alpha$ will be $\frac
  p{a_1}$ for some $1 \leq p<a_1$ with $(p,a_1)=1$ so there are $\varphi(a_1)$
  choices for $\rho_{\alpha}$. The angles of the rays $\mathcal{R}_{\alpha}$
  landing at $\alpha$ form a {\it rotation set} in the sense of \cite{Gold};
  that is, they form a finite subset of $\mathbb{R}/\mathbb{Z}$ that permutes
  cyclically under the circle map $\theta \mapsto d\theta (\text{mod }1)$. In
  \cite{Gold} it is shown that for given degree $d$ and rotation number
  $\rho_{\alpha}$ there are exactly $d-1$ disjoint rotation sets,
  distinguished by the relative position of their elements with respect to the
  $(d-1)^{\text{st}}$ roots of unity. Therefore, given the H-composition $P$
  with initial part $a_1$, there is a total of
  \begin{equation}
    \varphi(a_1)\cdot(d-1)
  \end{equation}
  choices for the set of angles of the rays $\mathcal{R}_{\alpha}$. By
  \cite{Portraits}, the widest angular gap between consecutive rays in
  $\mathcal{R}_{\alpha}$ corresponds to the 2 rays that delimit $Y_0$. Let us
  call these angles $\tau^-, \tau^+$.

  By the Douady-Hubbard theory, the ray $r_0$ with angle 0 lands at a fixed
  point $\beta \in Y_0$ different from $\alpha$. Choose a simple curve $\gamma
  \subset K_c$ joining the critical point 0 to $\beta$; then, the union of
  $\gamma$ and $r_0$ splits $Y_0$ in two parts. It is easy to see that the
  invariance relations (\ref{eqn:Foliations}), together with the normal form
  of $\phi_c$, force a well defined correspondence between the preimages of a
  ray $r_{\theta}$ under the $d$ branches of ${f_c^{-1}}_{|_{\mathbb{C}
  \setminus (\gamma \cup r_0)}}$, and the preimages of the angle $\theta$
  under the $d$ inverse branches of the circle map $\theta \mapsto d\theta\,
  (\text{mod }1)$.

  Each inverse branch of $\theta \mapsto d\theta\, (\text{mod }1)$ has the
  form $\frac{\theta+\kappa d}d$ with $0 \leq \kappa \leq d-1$. As a
  consequence, if we select $n-1$ consecutive branches of $f_c^{-1}$, the
  $(n-1)^{\text{st}}$ preimages of the angles $\tau^-$ and $\tau^+$ can be
  computed explicitly in terms of $\tau^-, \tau^+$:
  \begin{center}\begin{equation}\label{eqn:Etas}
    \eta^{\pm} :=
     \frac{\tau^{\pm} + \big( \kappa_1^{\pm} \cdot d \big) + \big(
     \kappa_2^{\pm} \cdot d^2 \big) + \ldots + \big( \kappa_{n-2}^{\pm} \cdot
     d^{n-2} \big)}{d^{n-1}} =
    \frac{\tau^{\pm} + X^{\pm}}{d^{n-1}},
  \end{equation}\end{center}
  where each coefficient $\kappa$ ranges between 0 and $d-1$ and is determined
  by the choice of branch. When a piece intersects the slit $\gamma \cup r_0$,
  its 2 rays are transformed by 2 different inverse branches of $\theta
  \mapsto d\theta\, (\text{mod }1)$. In particular, $\gamma \subset Y_0$
  implies that the first coefficients $\kappa_1^-$ and $\kappa_1^+$ are
  different, so $X^- \neq X^+$.

  It should be emphasized that the expressions to the right of $\tau_{\pm}$
  can be read as numbers written in base $d$. Then it is clear that every
  choice of inverse branches determines a different pair of angles $(\eta^-,
  \eta^+)$, because different itineraries encode different pairs $(X^-,\,
  X^+)$. \\

  Suppose $c$ is such that its itinerary is non-renormalizing (whether $P$ is
  renormalizable or not); then all pull-back pieces $Z_j$ are delimited by 2
  rays. In particular, $Z_1$ is delimited by the pair of rays $r_{\eta^-},\,
  r_{\eta^+}$. The particular sequence of inverse branches of $f_c^{-1}$ that
  is obtained by way of the pull-backs along $\text{Crit}_0$, is described
  next in terms of $P$.

  Recall that $\text{Crit}_0$ follows several circuits around the fixed point
  $\alpha$. Thus, among the components of $f_c^{-1}(Z_j)$, the unique
  candidate for $Z_{j-1}$ is the component that precedes $Z_j$ in the fan of
  pieces around $\alpha$. The only exception is, of course, when $Z_j$ is at
  the beginning of the current circuit around $\alpha$; i.e. when $Z_{j-1}$ is
  meant to be found somewhere in $Y_0$. In that case, $\zeta_k \in Z_{j-1}$
  for some $k$ and the number of candidate locations for $Z_{j-1}$ within
  $Y_0$ is either $d$ or $d-1$. The exact number of choices is determined by
  the value of $a_k$.

  Specifically, if $a_k=q$, then $Z_j \subset Y_1$ so all the $d$ components
  of $f_c^{-1}(Z_j)$ lie in $C \subset Y_0$ and satisfy the itinerary data,
  whereas if $a_k<q$ then $Z_j \subset Y_{q-a_{k+1}}$, and one component of
  $f_c^{-1}(Z_j)$ is in $Y_{q-a_k} \not\subset Y_0$. The other $d-1$
  components of $f_c^{-1}(Z_j)$ are located in $B_{(1,q-a_k)}, B_{(2,q-a_k)},
  \ldots, B_{(d-1,q-a_k)} \subset Y_0$ so the location of $Z_{j-1}$ (and the
  current choice of inverse branch of $f_c^{-1}$) can be encoded by the value
  $1 \leq b_k \leq d-1$.

  In the end, each value $a_k<q$ allows $d-1$ choices for $b_k$, translating
  into $d-1$ admissible inverse branches of $\theta \mapsto d\theta\,
  (\text{mod }1)$ at that step. If the itinerary is non-renormalizing, then
  each $a_k=q\, (k>1)$ results in a choice of $d$ possible branches since all
  components of $f_c^{-1}$ at that step lie in $C \subset Y_0$. The location
  of any other $Z_j$ is uniquely prescribed by $\rho_{\alpha}$.

  By definition, $\#\{a_k=q \mid k>1 \}=\omega(P)=w$, so we can write 
  $\#\{a_k<q\}=r-1-w$, where $r$ is the length of the H-composition $P$.
  The above discussion shows that
  \begin{gather}
    \#\big\{\text{itineraries associated to }
      P\big\}=(d-1)^{r-1-w}\, \text{ and}\label{eqn:Less_n_q} \\
    \#\big\{\text{pairs }(\eta^-,\,\eta^+)\text{ associated to a
     non-renormalizing itinerary}\big\}=d^w. \label{eqn:Equal_to_q}
  \end{gather}

  \noindent\underline{$P$ is non-renormalizable:} When the H-composition is
  not renormalizable, no itinerary can be renormalizing. It follows that if
  there is a $d$-center $c$ that satisfies $P(c)=P$, then the corresponding
  piece $Z_1$ is delimited by a pair of rays $(\eta^-,\, \eta^+)$ whose angles
  must belong to a collection of
  \begin{equation}
    \left[ \varphi(a_1)(d-1) \right] \cdot
    \left[ (d-1)^{r-1-w} \right] \cdot
    \left[ d^w \right]
  \end{equation}
  possible pairs. The first factor is given by an initial selection of rays
  to delimit $Y_0$, while the other 2 factors are consequence of
  (\ref{eqn:Less_n_q}), (\ref{eqn:Equal_to_q}) and the structure of $P$.

  It only remains to show that each admissible pair of angles $(\eta^-,\,
  \eta^+)$ will indeed contribute one $d$-center. \\

  Fix the pair $(\eta^-,\, \eta^+)$ described by a candidate sequence of
  pull-back choices. Recall that these choices include the selection of a pair
  $(\tau^-,\, \tau^+)$. By Formula (\ref{eqn:Etas}), there are 2 distinct
  linear functions that satisfy
  \begin{equation*}
    \tau^- \mapsto \frac{\eta^- + X^-}{d^{n-1}} \hspace{.15in}
    \text{and}  \hspace{.15in}
    \tau^+ \mapsto \frac{\eta^+ + X^+}{d^{n-1}}.
  \end{equation*}
  respectively. Let $(\theta^-,\,\theta^+)$ be the pair of fixed points of
  these functions:
  \begin{equation*}
    \frac{\theta^{\pm} + X^{\pm}}{d^{n-1}} = \theta^{\pm}.
  \end{equation*}
  Then, relative to the circle map $\theta \mapsto d\theta\, (\text{mod }1)$,
  the angles $\theta^-,\, \theta^+$ are periodic with period $n$, and have the
  required itinerary. Let $\Theta_0 = \{ \theta^-,\, \theta^+ \}$ and define
  recursively $\Theta_j$ as the images under $\theta \mapsto d\theta\,
  (\text{mod }1)$ of the 2 angles in $\Theta_{j-1}$ (for $j=1, \ldots, n-1$).
  By construction, the family $\big\{ \Theta_0, \ldots, \Theta_{n-1} \big\}$
  is a {\it formal orbit portrait} in the sense of \cite{Portraits}. That is,

  {\renewcommand{\theenumi}{{\bf \alph{enumi}}}
  \begin{enumerate}
    \item Each $\Theta_j$ is a finite subset of $\mathbb{R}/\mathbb{Z}$.

    \item For all $j \text{ mod } p$, $\theta \mapsto d\theta\,
          (\text{mod }1)$ maps $\Theta_{j-1}$ bijectively onto $\Theta_j$.

    \item All angles in $\Theta_0 \cup \Theta_1 \cup \ldots \cup \Theta_{n-1}$
          are periodic, with common period $kn$.

    \item The sets $\Theta_0,\, \Theta_1,\, \ldots,\, \Theta_{n-1}$ are
          pairwise unlinked; i.e. for $i \neq j$, $\Theta_i$ and $\Theta_j$
          are contained in disjoint intervals of $\mathbb{R}/\mathbb{Z}$ (by
          property \ref{Y-three}).
  \end{enumerate}
  }

  We have established that every non-renormalizing itinerary has associated a
  fixed number of formal orbit portraits. By Corollaries 5.4 and 5.5 of
  \cite{Portraits}, in the case $d=2$ there is a unique parameter $c_0$ with a
  parabolic periodic orbit $\mathcal{O}$ that follows the chosen itinerary,
  and such that each pair of rays $\Theta_j$ land at a common point of
  $\mathcal{O}$. Then, by Lemma 4.5 of \cite{Portraits}, there is a unique
  parameter $c$ associated to $c_0$ such that the critical orbit is periodic
  and has the chosen itinerary.

  This settles the case $d=2$; the description of orbit portraits for the
  general case $d \geq 2$ has never been explicitly written. However, as
  mentioned in the introduction of \cite{Portraits}, the case $d \geq 2$ does
  follow from assorted results in \cite{Eber} and \cite{Multibrots} that are
  equivalent to the treatment in \cite{Portraits}. Then the above discussion
  holds similarly for general $d$. \\

  \noindent\underline{$P$ is renormalizable:} In the case of a renormalizable
  H-composition, some itineraries are renormalizing and some are not. The
  argument presented above shows that a single non-renormalizing itinerary has
  associated
  \begin{equation*}
    \left[ \varphi(a_1)(d-1) \right] \cdot
    \left[ d^w \right]
  \end{equation*}
  $d$-centers. We will show that this count applies also for renormalizing
  itineraries. The reason for dealing with this case apart is that the piece
  $Z_0$ may contain points of $\text{Crit}_0$ other than 0. This would mean
  that the itinerary and related pull-back information are not enough to
  differentiate the behavior of similar critical orbits. This invalidates the
  previous counting argument. \\

  Recall that a renormalizing itinerary satisfies (\ref{eqn:Ren}) for some
  index $i<n$. Let us restrict attention to the smallest such index and call
  it $r'$. This has the effect that the H-subcomposition $P':=[a_1 + \ldots +
  a_{r'} = \frac{r'n}r]$ cannot be further decomposed in the manner of
  (\ref{eqn:Ren}). For notational consistency, define $w' := \omega(P')$
  (recall that $w:=\omega(P)$) and $n' := \frac{r'n}r$. Observe that $w =
  \frac{r(w'+1)}{r'} - 1$. \\

  Since $P'$ is not renormalizable, we already know that there are $\big[
  \varphi(a_1) \cdot d^{w'} \big]$ parameters satisfying $f^{(n')}(0) = 0$.
  Let us choose one and call it $c_0$. In \cite{p-l}, it is shown that there
  is a neighborhood $N$ of $c_0$, such that for every $c \in N$, $f_c$ follows
  the same itinerary as $f_{c_0}$ for $n'$ iterates and the polynomial
  $F_c(z):=f_c^{(n')}(z)$ maps $Z_0 \mapsto Z_{r'} \mapsto \ldots \mapsto
  Z_{n-r'} \mapsto Y_0$, in a $d^{w'+1}$ to 1 manner at every step. In fact,
  $F_c(z)$ is {\it polynomial-like} in the sense of \cite{p-l}; i.e. it maps
  the region\footnote{Strictly speaking, it is necessary that $Z_0$ is
  compactly contained in $Z_{r'}$: $Z_0 \Subset Z_{r'}$. This is not true in
  general, but $Z_0$ and $Z_{r'}$ can be slightly ''thickened'' to satisfy
  this stronger inclusion condition.} $Z_0 \subset Z_{r'}$ onto the larger
  region $Z_{r'}$ as a $(d^{w'+1})$-branched cover. In other words,
  ${F_c}_{|_{\mathbb{C}}}$ has global degree $d^{r'}$ but, when restricted to
  $Z_0$, it behaves like a degree $d^{w'+1}$ polynomial. The fundamental
  result of \cite{p-l} is \\

  \noindent
    {\bf The Straightening Theorem:} {\it Let $Z' \Subset Z$ be two open
    regions in $\mathbb{C}$ and $F:Z' \longrightarrow Z$ a polynomial-like map
    of degree $\delta$. Then $F$ is hybrid equivalent (quasi-conformally
    conjugate) to a polynomial $Q(z)$ of degree $\delta$. Moreover, if the
    filled Julia set of $Q$ is connected, then $Q$ is unique up to conjugation
    by an affine map.
  } \\

  Thus, there exists a degree $d^{w'+1}$ polynomial $Q_c$ associated to $F_c$,
  such that $Q_c$ has the same dynamics as ${F_c}_{|{Z_0}}$; i.e. the same
  behavior of critical orbits. Since ${F_c}_{|{Z_0}}$ has only one critical
  point, we can assume that the {\it straightening} $Q_c$ is the unique
  binomial $z^{(d^{w'+1})}+\hat{c}$ in the affine conjugacy class of $Q_c$.

  It follows from \cite{Orsay}-\cite{p-l} that $\{f_c \mid c \in N \}$ is a
  {\it full analytic family} of polynomial-like maps. As a consequence, any
  binomial $Q(z) = z^{(d^{w'+1})} + \hat{c}$ is the straightening of $F_c(z)$
  for a unique $c \in N$. \\

  Now, we are interested in parameters $c$ such that $f_c^n(0)=0$ or
  equivalently, $F_c^{(r/r')}(0)=0$. Clearly, any $(d^{w'+1})$-center
  $\hat{c}$ with period $\frac r{r'}$ represents the straightening $Q_c$ of
  some $F_c$ that satisfies $F_c^{(r/r')}(0)=0$. Moreover, since the critical
  point is periodic, the filled Julia set of $Q_c$ is connected. Then, by the
  Straightening Theorem the correspondence between $Q_c$ and $F_c$ is
  bijective.

  It only remains to count how many $d$-centers follow the given renormalizing
  itinerary. First, there is a choice among $\varphi(a_1)(d-1)$ rotation sets
  for the initial angles around $\alpha$. Next, since $f_c$ must follow the
  non-renormalizing H-subitinerary $P'$, there are $d^{w'}$ choices for the
  map $F_{c_0}$ (only those steps where $a_j=a_1$ admit choices since the
  value of the $b_j$'s is prescribed by $P'$). Finally, the neighborhood $N$
  of $c_0$ contains one parameter $c$ for every $\big( d^{w'+1} \big)$-center
  $\hat{c}$ with period $\frac r{r'}$; meaning that $f_c$ will satisfy
  $f_c^n(0) = f_c^{\big( \frac{n'r}{r'} \big)}(0) = F_c^{\big( \frac r{r'}
  \big)}(0) = 0$. By Lemma \ref{count} there are $\big( d^{w'+1} \big)^{\frac
  r{r'}+1}$ such parameters $\hat{c}$; since $(w'+1)(\frac r{r'}-1) =
  \frac{r(w'+1)}{r'} - 1 - w' = w - w'$, the total of $d$-centers associated
  to the renormalizing itinerary is
  \begin{gather*}
    \left[ \varphi(a_1)(d-1) \right]
    \cdot\left[ d^{w'} \right]
    \cdot\left[ (d^{w'+1})^{\frac r{r'}-1} \right] = \\
    \left[ \varphi(a_1)(d-1) \right]
    \cdot\left[ d^w \right].
  \end{gather*}
\QED \\

Besides providing an intuitive frame for Formula (\ref{eqn:Formula}), the
above proof emphasizes the difference between renormalizable and
non-renormalizable H-compositions. The Formula should reflect interesting
number-theoretical properties related to the renormalization phenomenon; a
feature that is not immediately apparent in the proof of Section
\ref{sec:pf2}. This line of investigation, along with the search for similar
identities in other classes of parameters are the subject of forthcoming
research announcements by the second author.

\section{A Combinatorial Proof}\label{sec:pf2}

We give now a self-contained proof of Theorem \ref{Thm:Formula} based on the
more usual tools of partition theory. First we establish three elementary
identities.

\begin{lemma}\label{lemma} Let $\mathcal{C}(n,b,s)$ denote the number of
  compositions of $n$ into $b$ parts, each less than or equal to $s$. Then

\begin{enumerate}
  \item The following formal power series identity holds:
        \begin{equation*}
          \sum_{m=1}^{\infty} \frac{\varphi(m)z^m}{1-z^m} = \frac z{(1-z)^2}.
        \end{equation*}
  \item The generating function of $\mathcal{C}(n,b,s)$ with index $n$ is
        \begin{equation*}
          \sum_{n=0}^{\infty} \mathcal{C}(n,b,s)z^n = (1-z)^{-b}z^b(1-z^s)^b.
        \end{equation*}
  \item The number of H-compositions $P$ with first part $m$, multiplicity
        $\omega(P) = w$ and with a total of $r$ parts is
        \begin{equation*}
          {r-1 \choose w} \mathcal{C}\big( n-(w+1)m,r-w-1,m-1 \big).
        \end{equation*}
\end{enumerate}
\end{lemma}

\pf{s}
\begin{enumerate}
  \item This is a well known result on Lambert series; see Theorem 309 of
    \cite{Hardy}. It is proved by expanding $\frac1{1-z^m}$ as a formal
    geometric series:
    \begin{equation*}
      \sum_{m=1}^{\infty} \frac{\varphi(m)z^m}{1-z^m} =
      \sum_{m=1}^{\infty} \sum_{n=1}^{\infty} \varphi(m)z^{mn}.
    \end{equation*}
    Since $m$ divides $mn,$ this double sum can be rearranged as
    \begin{equation*}
      \sum_{N=1}^{\infty} \Big( \sum_{m \mid N} \varphi(m) \Big) z^N =
      \sum_{N=1}^{\infty} N z^N = \frac z{(1-z)^2}.
    \end{equation*}
  \item To each composition of $n$ in $b$ parts, $n=a_1+ \ldots +a_b$ with
    every part $a_j \leq s,$ corresponds a monomial $z^n = z^{a_1} \ldots
    z^{a_b}$ in the expansion of the product $(z+z^2+ \ldots + z^s)^b.$The
    coefficient of $z^n$ in $(z+z^2+ \ldots + z^s)^b$ is clearly
    $\mathcal{C}(n,b,s),$ so
    \begin{equation*}
      \sum_{n = 0}^{\infty} \mathcal{C}(n,b,s)z^n =
      (z+z^2+ \ldots + z^s)^b =
      \left( \frac{z-z^{s+1}}{1-z} \right)^b =
      (1-z)^{-b}z^b(1-z^s)^b.
    \end{equation*}
  \item Since $\omega(P) = w,$ there are $w+1$ parts $a_1,a_{s_1}, \ldots,
    a_{s_w}$ equal to $m$. The remaining parts are all smaller or equal than
    $m-1$ and form a composition of $n-(w+1)m$. Naturally, there are ${r-1
    \choose w}$ possible places to put the $a_{s_j}$. Hence the total number
    of prescribed compositions is ${r-1 \choose w} \mathcal{C}\big(
    n-(w+1)m,r-w-1,m-1 \big)$ as stated. \QED
\end{enumerate}

Backed by these results, we can proceed to prove the identity. \\

\pf{ of Theorem \ref{Thm:Formula}}
  We will show that the generating function
  \begin{equation*}
    G(z) = \sum_{n=1}^{\infty}
             \left( \sum_{\substack{P \in \mathcal{H}(n) \\ a_1+\ldots+a_r=n}}
               \varphi(a_1) \cdot (d-1)^{r-\omega(P)} \cdot d^{\omega(P)}
             \right) z^n
  \end{equation*}
  coincides with the power series
  $\sum_{n=1}^{\infty} (d^n-1)z^n = \frac{(d-1)z}{(1-d z)(1-z)}.$
  Then, a term by term comparison of the coefficients in both series yields
  the result. Start by representing $G(z)$ as follows:
  \begin{equation*}
    \sum_{n=1}^{\infty}
    \left(
      \sum_{m=1}^n \sum_{\substack{w \geq 0 \text{ s.t.}\\
                                   (w+1)m \leq n}}
      \sum_{r \geq w+1}
        {r-1 \choose w}\mathcal{C}\big( n-(w+1)m,r-w-1,m-1 \big) \cdot
      \varphi(m) (d-1)^{r-w} d^w
    \right) z^n
  \end{equation*}
  where the index $m$ represents all possible values of $a_1$ in a
  H-composition of $n;$ $w$ runs over the possible multiplicities of such
  compositions and $r$ stands for the lengths of compositions with the given
  $m$ and $w.$ Note that for every choice of parameters we count a total of
  ${r-1 \choose w}\mathcal{C}\big( n-(w+1)m,r-w-1,m-1 \big)$ compositions
  according to point 3 of Lemma \ref{lemma}.

  Replace the innermost index $r$ with $r+w+1.$ Also, observe that $m>n$ or
  $w>\frac nm - 1$ allow no valid H-compositions, so we are free to let the
  indices $m$ and $w$ run to $\infty$ and interchange the summation order:
  \begin{equation*}
    \sum_{m=1}^{\infty} \sum_{w=0}^{\infty} \sum_{n=1}^{\infty}
    \sum_{r=0}^{\infty} {r+w \choose w}\mathcal{C}\big( n-(w+1)m,r,m-1 \big)
      \cdot \varphi(m)(d-1)^{r+1}d^w \cdot z^m.
  \end{equation*}

  A similar simplification occurs when we replace $n$ by $n+(w+1)m.$ Since
  there are no compositions of negative numbers, we can let the sum over $n$
  start at $n=0:$
  \begin{equation*}
    \sum_{m=1}^{\infty} \sum_{w = 0}^{\infty} \sum_{n = 0}^{\infty}
    \sum_{r = 0}^{\infty} {r+w \choose w}\mathcal{C}\big( n,r,m-1 \big) \cdot
      \varphi(m) (d-1)^{r+1} d^w \cdot z^{n+(w+1)m}.
  \end{equation*}

  Now we proceed to eliminate the interior sums. First, by point 2 of Lemma
  \ref{lemma} we get:
  \begin{equation*}
    \sum_{m=1}^{\infty} \sum_{w = 0}^{\infty} \sum_{r = 0}^{\infty} 
    {r+w \choose w}\left( (1-z)^{-r}z^r(1-z^{m-1})^r \right) \cdot
    \varphi(m)(d-1)^{r+1}d^w \cdot z^{(w+1)m}.
  \end{equation*}

  Then, since $\sum {r+w \choose w} q^w = (1-q)^{-r-1}\,,$
  \begin{gather*}
    \sum_{m=1}^{\infty} \sum_{r = 0}^{\infty}
    \left( (1-z)^{-r}z^r(1-z^{m-1})^r \right) \cdot
    \varphi(m) (d-1)^{r+1} z^m (1-d
     z^m)^{-r-1}. \\
    \intertext{Gathering powers of $r$ together:}
    \sum_{m=1}^{\infty} \sum_{r = 0}^{\infty}
     \left( \frac{(d-1)z(1-z^{m-1})}{(1-d z^m)(1-z)} \right)^r \cdot
     \frac{\varphi(m) (d-1) z^m}{(1-d z^m)} = \\
    \sum_{m=1}^{\infty}
     \frac1{1-\frac{(d-1)z(1-z^{m-1})}{(1-d z^m)(1-z)}} \cdot
     \frac{\varphi(m) (d-1) z^m}{(1-d z^m)} = \\
    \sum_{m=1}^{\infty}
     \frac{\varphi(m)(d-1)z^m}{(1-dz^m)-\frac{(d-1)z(1-z^{m-1})}{(1-z)}} = \\
    (d-1)(1-z) \sum_{m=1}^{\infty}
     \frac{\varphi(m)z^m}{(1-dz^m)(1-z)-(d-1)z(1-z^{m-1})}.
    \intertext{The denominator simplifies to $1-z-dz^m+dz^{m+1}-dz+dz^m+z-z^m
     = (1-dz)(1-z^m),$ so the expression becomes:}
    \frac{(d-1)(1-z)}{(1-dz)} \sum_{m=1}^{\infty}
     \frac{\varphi(m) z^m}{1-z^m}.
    \intertext{Finally, the first point of Lemma \ref{lemma} gives:}
    G(z) =
     \frac{(d-1)(1-z)}{(1-dz)}
     \left( \frac z{(1-z)^2} \right) = \\
    \frac{(d-1)z}{(1-dz)(1-z)}.
  \end{gather*}
\QED


\begin{thebibliography}{99999}

\bibitem[A]{Andrews} G. E. Andrews,
  The Theory of Partitions.
  (Cambridge University Press, 1984).

\bibitem[D]{Dedekind} R. Dedekind,
  Theory of Algebraic Integers.
  (Cambridge University Press, 1996).

\bibitem[DH1]{Orsay} A. Douady \& J. H. Hubbard,
  \'Etude dynamique des polyn\^omes complexes I \& II.
  Publ. Math. Orsay, 1984-85.

\bibitem[DH2]{p-l} A. Douady \& J. H. Hubbard,
  On the dynamics of polynomial-like maps.
  Ann. Sci. \'Ec. Norm. Sup., 18 (1985), 287-343.

\bibitem[E]{Eber} D. Eberlein,
  Rational parameter rays of Multibrot sets.
  Diploma Thesis. Technische Universit\"at M\"unchen. \\
  Available at: {\tt ftp://ftp.math.sunysb.edu/theses/thesis99-2/part1.ps.gz}

\bibitem[G]{Gold} L. Goldberg,
  Fixed points of polynomial maps, I.
  Ann. Sci. \'Ec. Norm. Sup., 25 (1992), 679-685.

\bibitem[GM]{Gold2} L. Goldberg \& J. Milnor,
  Fixed points of polynomial maps, II.
  Ann. Sci. \'Ec. Norm. Sup., 26 (1992), 51-98.

\bibitem[HW]{Hardy} G. H. Hardy \& E. M. Wright
  An introduction to the theory of numbers.
  (Oxford University Press, 1979, $5^{\text{th}}$ edition).

\bibitem[HY]{Tableaux} J. H. Hubbard,
  Local connectivity of Julia sets and bifurcation loci: Three theorems of
  J.-C. Yoccoz.
  In: Topological Methods in Modern Mathematics pp. 467-511 (ed. L. Goldberg
  \& A. Phillips), (Publish or Perish, 1993).

\bibitem[LS]{Dierk} E. Lau \& D. Schleicher,
  Internal addresses in the Mandelbrot set and irreducibility of
  polynomials.
  Preprint \# 1994/19 IMS, Stony Brook.

\bibitem[M1]{M_book} J. W. Milnor,
  Dynamics in One Complex Variable.
  (Vieweg, 1999).

\bibitem[M2]{Portraits} J. W. Milnor,
  Periodic orbits, external rays and the Mandelbrot set: An expository
  account.
  In: Asterisque 261 `Geometrie Complexe et Systemes Dynamiques', pp.
  277-333, (SMF 2000).

\bibitem[P]{Poir} A. Poirier,
  On the realization of Fixed point portraits. (An addendum to ``Fixed point
  portraits'' by Goldberg and Milnor).
  Preprint \# 1991/20 IMS, Stony Brook.

\bibitem[S]{Multibrots} D. Schleicher,
  On fibers and local connectivity of Mandelbrot and Multibrot sets.
  Preprint \# 1998/13(a) IMS, Stony Brook.

\end{thebibliography}
\end{document}